\documentclass[11pt]{article}

\usepackage{inputenc}

\UseRawInputEncoding

\usepackage{latexsym,amssymb,amsmath,amscd,amsthm,amsxtra}
\usepackage{mathrsfs}
\usepackage{epsfig}
\usepackage{pgf,tikz}
\usepackage{graphicx}
\usepackage{graphics}
\usepackage{setspace}
\newcommand{\D}{\displaystyle}
\newcommand{\counte}{section}
\newtheorem{define}{\bf Definition}[\counte]

\newtheorem{lemma}{\bf Lemma}[\counte]
\newtheorem{theorem}{\bf Theorem}[\counte]
\newtheorem{coro}{\bf Corollary}[\counte]
\newtheorem{example}{\bf Example}[\counte]

\author{Zhao Xu-an, zhaoxa@bnu.edu.cn\\Gao Hongzhu, hzgao@bnu.edu.cn\\
Department of Mathematics, Beijing Normal University\\Key Laboratory
of Mathematics and Complex Systems\\ Ministry of Education,
China, Beijing 100875}

\title{General affine differential geometry of surfaces in affine space $A^3$, I: the elliptical case\thanks{The authors are supported by NSFC 12071034}}
\date{}

\begin{document}

\maketitle
\begin{abstract}
In this paper we study the general affine differential geometry of surfaces in affine space $A^3$. For a regular elliptical
surface we design an algorithm to compute the moving frame of minimal order and get the complete system of differential invariants.
As an application we classify regular elliptical surfaces of constant curvature up to affine congruence. The work in this paper shows that the essential ground of a Kleinian differential geometry is not analysis but algebra and hence formal symbolic computation play an important role here.
\end{abstract}

% All of the discussion can be implemented as an algorithm which can be applied to other Kleinian differential geometry.


\section{Introduction}

The affine geometry was founded by Blaschke, Pick, Radon, Berwald and Thomsen among others in the period
from 1916 to 1923. For accounts and expository books
appeared on the subject, see Blaschke\cite{Blaschke_85}, Guggenheimer\cite{Guggenheimer_63} and Spivak's\cite{Spivak}. In these work affine geometry means the equi-affine geometry. That is the Kleinian geometry of the affine transformation group which preserves volume. But for general affine geometry, i.e. the Kleinian geometry of general affine transformation group, there is little work. All the work we can find is Weise\cite{Weise_38}\cite{Weise_39}, Kllingenberg\cite{Klingenberg_81_1}\cite{Klingenberg_51_2}, Svec\cite{Svec_59},
Wilkinson\cite{Wilkinson_88} and Weiner\cite{Weiner_94}. The property of volume preserving makes things easier to be dealt with. Compared with the equi-affine geometry, the study of general affine geometry is more difficult.

%We need a set theory of general affine geometry which is parallel to the classical Euclidean geometry of curves and surfaces.

In Zhao and Gao\cite{Zhao_Gao161} we constructed the theory of general affine differential geometry of plane curves. In this paper we will consider the case for space surfaces. In the study of classical differential geometry,
it seems that the geometry and analysis are vital. We will show in this paper that neither geometry nor analysis but algebra is essential in the theory of general Kleinian differential geometry. This means that we can construct the Kleinian differential geometry theory by using pure symbolic computation methods. The only ingredient from analysis is that we assume the equations of surfaces we considered can be expanded as power series locally. The viewpoint in this paper makes the structure of a Kleinian differential geometry theory more transparent.

The input of our method is the action of space affine transformation group $Aff(3)$ on three dimensional affine space $A^3$.
The procedure to construct affine differential theory is as follows:

1. We determine the action of the affine transformation group on jet spaces of surfaces in $A^3$. Then the differential invariants of surfaces are converted to algebraic invariants of jets.

2. We choose the standard forms of jets of surfaces. By Fels and Olver\cite{Fels_Olver_98}\cite{Fels_Olver_99}\cite{Fels_Olver_04}, we have an equivariant moving frame.

3. Based on steps 1 and 2, by using Cartan's moving frame method we construct the affine differential geometry theory of elliptical surfaces.

This strategy is valid for any Kleinian differential geometry.

In this paper we need two algorithms to finish the computation, we call them algorithm A1 and A2.

A1. An algorithm to determine the standard form of a jet of a surface. See section 3.

%2. An algorithm to expand the infinitesimal generators(or vector fields) of $G$ act on $A^3$ to the infinitesimal generators(or vector fields) of $G$ act on jet space %$J^{2,3}(A^3)$. See section 4.

A2. An algorithm to compute the moving equations of the moving frame and compatible conditions.
%The input data is the standard form of a jet and expanded vector fields on jet space $J^{2,3}(A^3)$. See Theorem 4.1 in section 4.

This paper is part of a project to study the algebraic structure of general Kleinian differential geometry theory. For simplicity, we only consider the regular elliptical surfaces here. The hyperbolic and parabolic cases will
 be considered elsewhere.

 The contents of this paper are as follows: In section 2, we give an introduction to jet spaces and Fels and Olver's equivariant moving frame method. In section 3, we design an algorithm to compute the standard forms
 of jets of surfaces under the space affine transformation group.
 In section 4, we construct the differential geometry theory of regular elliptical surfaces. We give an algorithm to compute the moving equations of equivariant moving frames and to compute the compatible conditions. As a result we obtain the complete system of
 differential invariants. In section 5, as an example, we give an explicit surface to compute its moving frame and moving equation using our algorithm. We also classify elliptical surfaces with constant curvature up to
 affine congruence. Section 6 is an appendix which gives the details of solving the compatible equations.
 %And We define the second, third and fourth fundamental form for elliptical space surfaces and show that they determine the congruence class of a surface.
%In section 6 we give an affine classification curves of constant curvature. In section 7 we discuss modular invariants.certain affine invariant differential forms on a surface and give construct

\section{Jet spaces and the equivariant moving frame method}
This section is an introduction to jet spaces and Fels and Olver's equivariant moving frame method. See Fels and Olver\cite{Fels_Olver_98}\cite{Fels_Olver_99}\cite{Fels_Olver_04} for details.

\subsection{Jet spaces}
The differential geometry of a submanifold $N$ at a point $P$ is determined by the shape of $N$ at an arbitrary small neighborhood. So all the local differential geometric properties and invariants of $N$ are determined by the local
data of $N$. It is useful to isolate the information of $N$ at $P$. This idea hints the concept of jet of submanifold. The general definition of jet of submanifolds was given by Ehresmann in 1950s for an ambient manifold $M$.
%See Ehresmann\cite{Olver_95} for reference.

For a smooth manifold $M$ of dimension $n$, let $SM^{d}_P(M)$ be the set of all smooth $d$-dimensional submanifolds of $M$ that contain the point $P$ and $SM^{d}(M)=\bigcup\limits_{P\in M} SM^{d}_P(M)$.
For integer $r\geq 0$, define an equivalence relation $\sim$ on $SM^{d}_P(M)$ such that $N_1\sim N_2$ iff $N_1$ and $N_2$ have contact of order $r$. The jet space of $d$-dimensional submanifolds of $M$ at a point $P$ of order $r$ is the quotient set $J_P^{d,r}(M)=SM^{d}_P(M)/\sim$. And $J^{d,r}(M)=\bigcup\limits_{P\in M} J_P^{d,r}(M)$ is the jet space of $d$-dimensional submanifolds of $M$ of order $r$. An element of jet space $J^{d,r}(M)$ is called a $d$-jet of order $r$. For an $n$-dimensional manifold $M$, $J^{d,r}(M)$ is a manifold of dimension $\displaystyle{d+(n-d){d+r\choose r}}$. In this paper we use the Jet space $J^{2,r}(A^3)$ for $r>0$ in the study of surfaces in affine space $A^3$.

Jet space can be regarded as the finite dimension cut-off of infinite dimension space of all submanifolds with fixed dimension. The use of jet spaces separates the study of local differential geometry into algebra part and analysis part. So it makes the structure of local differential geometry theory more transparent.
\subsection{Fels and Olver's moving frame method}

\begin{define} Let $\phi:G\times M\to M$ be a smooth action of Lie group $G$ on smooth manifold $M$. A moving frame on $M$ is a smooth, G-equivariant map $\rho: M\to G$.
\end{define}

There are two types of moving frames.

$\left\{
  \begin{array}{ll}
    \rho(gz) =g\rho(z), & \hbox{left moving frame;} \\
    \rho(gz) =\rho(z)g^{-1}, & \hbox{right moving frame.}
  \end{array}
\right.$

\begin{theorem}(Fels and Olver)
A moving frame exists on $M$ if and
only if $G$ acts freely and regularly on $M$.
\end{theorem}

The explicit construction of a moving frame is based on Cartan's normalization procedure.
Let $G$ act freely and regularly on $M$ and $K$ be a cross-section to the group orbits, that is a submanifold $K$ which transversally intersects each orbit once. Let $g$ be the unique group element which maps $P$ into the
cross-section $K$, then $\rho:M \to G,P\mapsto g$ is a right moving frame. And $\rho: M \to G,P\mapsto g^{-1}$ is a left moving frame.
The unique intersection point of the orbit of $P$ and $K$ can be regarded as the standard form of $P$, as
prescribed by the cross-section $K$.

If a moving frame is in hand, the determination of the invariants is routine.
The specification of a moving frame by choosing a cross-section induces a canonical
procedure to map functions on $K$ to invariants.

\begin{define} The invariantization of a function $F: M\to\mathbb{R}$ is the unique invariant
function $\iota(F)$ that coincides with $F$ on the cross-section, that is $\iota(F)|_K= F|_K$.\end{define}

Invariantization defines a projection from the space of smooth functions to the space of invariants that, moreover,
preserves all algebraic operations. The fundamental differential invariants are obtained by invariantization of coordinate functions on jet space.

%Since the action $Aff(2)\times J^{1,r}(A^2)\to J^{1,r}(A^2)$ is not free, so the moving frame can only be constructed on certain open set of $J^{1,r}(A^2)$. We will discuss this later.

\section{The action of $Aff(3)$ on $J^{2,r}(A^3)$}
Let $Aff(3)$ be the general affine transformation group of affine space $A^3$. %The dimension of $Aff(2)$ is $6$.
The general affine geometry is the Kleinian geometry given by the natural group action $Aff(3)\times A^3\to A^3$. Under the affine coordinates $x,y,z$ on $A^3$, a general affine transformation $T$ has the form
$$\left(
\begin{array}{c}
 x \\
  y \\ z\\
   \end{array}
   \right)=\left(
             \begin{array}{ccc}
               b_{11} & b_{12} & b_{13}\\
               b_{21} & b_{22} & b_{23}\\
               b_{31} & b_{32} & b_{33}\\
             \end{array}
           \right)
   \left(
\begin{array}{c}
 x' \\
  y' \\ z'\\
   \end{array}
   \right)+\left(
\begin{array}{c}
 x_0 \\
  y_0\\
  z_0\\
   \end{array}
   \right).$$
   Or
$$\left(
\begin{array}{c}
 x \\
  y \\
  z\\
  1 \\
   \end{array}
   \right)=\left(
             \begin{array}{cccc}
               b_{11} & b_{12} & b_{13} & x_0\\
               b_{21} & b_{22} & b_{23} & y_0\\
               b_{31} & b_{32} & b_{33} & z_0\\
               0&0&0&1\\
             \end{array}
           \right)
   \left(
\begin{array}{c}
 x' \\
  y' \\
  z'\\
   1\\
   \end{array}
   \right).$$

\noindent Where $\det \left(
             \begin{array}{ccc}
               b_{11} & b_{12} & b_{13} \\
               b_{21} & b_{22} & b_{23} \\
               b_{31} & b_{32} & b_{33} \\
                            \end{array}
           \right)\not=0$. Hence the coordinates on $Aff(3)$ are $b_{11},b_{12},b_{13},b_{21},b_{22},b_{23},b_{31},$
           $b_{32},b_{33}$, $x_0,y_0,z_0$.

In the following we assume that a surface has the form $z=z(x,y)$ locally. The local coordinates of the jet space $J^{2,3}(A^3)$ are $x,y,z,z_x,z_y,z_{xx},z_{xy},z_{yy}, z_{xxx},z_{xxy},z_{xyy} ,z_{xyy},z_{yyy}$. The action of $Aff(3)$ on $A^3$ induces an action on $SM^{2}(A^3)$ and hence on $J^{2,3}(A^3)$.

We expand $z(x,y)$ at $(x_0,y_0)$ as
$$z(x,y)=a_{00}+a_{10}(x-x_0)+a_{01}(y-y_0)+\frac{1}{2}(a_{20}(x-x_0)^2+2a_{11}(x-x_0)(y-y_0)+a_{02}(y-y_0)^2)$$
\begin{equation}
+\frac{1}{6}(a_{30}(x-x_0)^3+3a_{21}(x-x_0)^2(y-y_0)+3a_{12}(x-x_0)(y-y_0)^2+a_{03}(y-y_0)^3)+\cdots.
\end{equation}
%$$+\frac{1}{6}(a_{40}x^4+4a_{31}x^3y+6a_{22}x^2y^2+4a_{13}xy^3+a_{04}y^4).$$%Since any jet can be transformed by elements in $Aff(2)$ into this form, it gives no essential restriction.

We identify the symbols $z,z_x,z_y,z_{xx},z_{xy},z_{yy}, z_{xxx},z_{xxy},z_{xyy} ,z_{xyy},z_{yyy},\cdots $ with $a_{00}, a_{10},a_{01},a_{20},$ $a_{11},a_{02},a_{30},a_{21},a_{12},a_{03},\cdots$.

In the following, we will use affine transformations to transform a jet of surface $S:z(x,y)$ at $P(x_0,y_0,z_0)$ to a jet of surface $S': z'(x',y')$ at $P'(x'_0,y'_0,z'_0)$ step by step, such that the final jet of surface has certain standard form. To simplify our notations, we often use the same symbols $x,y,z$ instead of $x',y',z'$ after the affine transformations.

First we observe that a jet as in Equation 1 can be transformed by a translation to a jet at $(0,0,0)$ which has the form
%$$z=u'(x,y)$$
\begin{equation} z(x,y)=a_{10} x+a_{01} y+\frac{1}{2}(a_{20}x^2+2a_{11}xy+a_{02}y^2)+\cdots.\end{equation}

%To simplify the notation, we rewrite the symbol $x',y',z'$ by $x,y,z$ in Equation 3.

Let $T_1$ be an affine transformation $x=b_{11}x'+b_{12}y'+b_{13}z', y=b_{21}x'+b_{22}y'+b_{33}z', z=b_{31}x'+b_{32}y'+b_{33}z'$ preserving the point $O(0,0,0)$.
\begin{lemma}
Under the affine transformation $T_1$, if a jet of the form in Equation 2 is transformed to the form
%\begin{equation}z(x,y)=a_{10} x+a_{01} y+\frac{1}{2}(a_{20}x^2+2a_{11}xy+a_{02}y^2)+\cdots\end{equation}

$$z'(x',y')=a'_{10} x'+a'_{01} y'+\frac{1}{2}(a'_{20}x'^2+2a'_{11}x'y'+a'_{02}y'^2)+\cdots,$$
\noindent then we have
\begin{equation}(b_{33}-a_{10}b_{13}-a_{01}b_{23})(a'_{10}, a'_{01})= (b_{31}, b_{32})-(a_{10}, a_{01})\left( \begin{array}{cc}
                     b_{11} & b_{12} \\
                     b_{21} & b_{22} \\
 \end{array}
 \right).\end{equation}

%\begin{equation}
%a'_{10}=-\frac{b_{31}-a_{10}b_{11}-a_{01}b_{21}}{b_{33}-a_{10}b_{13}-a_{01}b_{23}}\label{a10}
%\end{equation}

%\begin{equation}
%a'_{01}=-\frac{b_{32}-a_{10}b_{12}-a_{01}b_{22}}{b_{33}-a_{10}b_{13}-a_{01}b_{23}}\label{a01}
%\end{equation}
\end{lemma}

\noindent{\bf Proof}: This is a direct computation from the affine transformation $T_1$ and Equation 2.
\begin{lemma}
(1)We can choose affine transformations $T_1$ suitably such that a jet of form in Equation 2 is transformed to a form with $a_{00}=a_{10}=a_{01}=0$

(2)The subgroup $G_1$ of $Aff(3)$ keeping the condition $a_{00}=a_{10}=a_{01}=0$ is

 $\D{G_1=\left\{ \left(
                   \begin{array}{ccc}
                     b_{11} & b_{12} & b_{13} \\
                     b_{21} & b_{22} & b_{23} \\
                     0 & 0 & b_{33} \\
                   \end{array}
                 \right)
\mid(b_{11}b_{22}-b_{12}b_{21})b_{33}\not=0, \right\}}$.
\end{lemma}

\noindent{\bf Proof}: (1) It is easy to see that we can choose $b_{11},b_{21},b_{31},b_{12},b_{22},b_{32}$ with $b_{31}-a_{10}b_{11}-a_{01}b_{21}=0$ and $b_{32}-a_{10}b_{12}-a_{01}b_{22}=0$. Then by Equation 3 we have $a'_{10}=a'_{01}=0$.

(2) This can be seen from Equation 3 as we have now $a_{00}=a_{10}=a_{01}=a'_{00}=a'_{10}=a'_{01}=0$.

%This is a direct result of Equation \ref{a10} and Equation \ref{a01}.
%Now we rewrite the symbol $x',y',z'$ by $x,y,z$.

Let $T_2$ be an affine transformation in the subgroup $G_1$ of form $x=b_{11}x'+b_{12}y'+b_{13}z', y=b_{21}x'+b_{22}y'+b_{33}z', z=b_{33}z'$.

\begin{lemma}
Under the affine transformation $T_2$, if a jet of the form \begin{equation}z(x,y)=\frac{1}{2}(a_{20}x^2+2a_{11}xy+a_{02}y^2)+\cdots \end{equation}
is transformed to the form
$$z'(x',y')=\frac{1}{2}(a'_{20}x'^2+2a'_{11}x'y'+a'_{02}y'^2)+\cdots,$$
\noindent then we have

\begin{equation}b_{33}\left( \begin{array}{cc}
                     a'_{20} & a'_{11} \\
                     a'_{11} & a'_{02} \\
 \end{array} \right)
 =\left( \begin{array}{cc}
                     b_{11} & b_{21} \\
                     b_{12} & b_{22} \\
 \end{array} \right)
 \left( \begin{array}{cc}
                     a_{20} & a_{11} \\
                     a_{11} & a_{02} \\
 \end{array} \right)
 \left( \begin{array}{cc}
                     b_{11} & b_{12} \\
                     b_{21} & b_{22} \\
 \end{array}
 \right)
\end{equation}

%\begin{equation}
%b_{33}a'_{20}=a_{20}b_{11}^2+2a_{11}b_{11}b_{21}+a_{02}b_{21}^2
%\end{equation}
%\begin{equation}
%b_{33}a'_{11}=a_{20}b_{11}b_{12}+a_{11}b_{11}b_{22}+a_{11}b_{12}b_{21}+a_{02}b_{21}b_{22}
%\end{equation}
%\begin{equation}
%b_{33}a'_{02}=a_{20}b_{12}^2 +2a_{11}b_{12}b_{22}+a_{02}b_{22}^2
%\end{equation}
\end{lemma}

\noindent{\bf Proof:} This follows directly from $T_2$ and Equation 4.

%In the following we divide the discussion into two cases.

%\subsection{The elliptical case}
By rotating in the $xy$-plane and re-scaling on axes $x$ and $y$, the quadratic items of a jet can be transformed into one of the four standard forms $\frac{1}{2}(x^2+y^2),\frac{1}{2}(x^2-y^2),\frac{1}{2}x^2$ or $0$.
\begin{define}We call the corresponding four jet types the elliptical, hyperbolic, parabolic and degenerate types.\end{define}
In jet spaces $J^{2,3}(A^3)$,
all jets of elliptical type form an open sub-manifold. We call a point $P$ in a surface $S$ is elliptical if the jet of surface $S$ at $P$ is elliptical.
In this paper, we only consider elliptical surfaces, i.e. surfaces whose points are all elliptical. %The hyperbolic surface can be dealt with similarly.

\begin{coro}
The jet of elliptical surface $S$ can be transformed to the standard form of $$z(x,y)=\frac{1}{2}(x^2+y^2)+\frac{1}{6}(a_{30}x^3+3a_{21}x^2 y+3a_{12}xy^2+a_{03}y^3)+\cdots$$
\end{coro}

\begin{lemma}
Let $G_2$ be the subgroup of $G_1$ keeping the condition $a_{00}=a_{10}=a_{01}=a_{11}=0,a_{20}=a_{02}=1$, then $\D{G_2=\left\{ \left(
                   \begin{array}{ccc}
                     b_{11} & b_{12} & b_{13} \\
                     b_{21} & b_{22} & b_{23} \\
                     0 & 0 & b_{33} \\
                   \end{array}
                 \right)
\mid b^2_{11}+b^2_{21}=b^2_{12}+b^2_{22}=b_{33},b_{11}b_{12}+b_{21}b_{22}=0 \right\}}$.\label{G2}
\end{lemma}

\noindent{\bf Proof:} This follows from Equation 5. In the present case $a_{11}=a'_{11}=0, a_{20}=a_{02}=a'_{20}=a'_{02}=1$.

Let $T_3$ be an affine transformation $x=b_{11}x'+b_{12}y'+b_{13}z', y=b_{21}x'+b_{22}y'+b_{33}z', z=b_{33}z'$ which is in $G_2$,

\begin{lemma}
Under the affine transformation $T_3$, if a jet of the form \begin{equation}z(x,y)=\frac{1}{2}(x^2+y^2)+\frac{1}{6}(a_{30}x^3+3a_{21}x^2 y+3a_{12}xy^2+a_{03}y^3)+\cdots\end{equation} \noindent is transformed to the form
$$z'(x',y')=\frac{1}{2}(x'^2+y'^2)+\frac{1}{6}(a'_{30}x'^3+3a'_{21}x'^2 y'+3a'_{12}x'y'^2+a'_{03}y'^3)+\cdots, $$ \noindent then we have
\begin{equation}
b_{33}a'_{30}=a_{30}b_{11}^3+3a_{21}b_{11}^2b_{21}+3a_{12}b_{11}b_{21}^2+a_{03}b_{21}^3+3b_{11}b_{13}+3b_{21}b_{23}
\end{equation}
\begin{equation}
b_{33}a'_{21}=a_{30}b_{11}^2b_{12}+a_{21}b_{11}(b_{11}b_{22}+2b_{12}b_{21})+a_{12}b_{21}(2b_{11}b_{22}+b_{12}b_{21})+a_{03}b_{21}^2b_{22}+b_{12}b_{13}+b_{22}b_{23}
\end{equation}
\begin{equation}
b_{33}a'_{12}=a_{30}b_{11}b_{12}^2+a_{21}b_{12}(2b_{11}b_{22}+b_{12}b_{21})+a_{12}b_{22}(b_{11}b_{22}+2b_{12}b_{21})+a_{03}b_{21}b_{22}^2+b_{11}b_{13}+b_{21}b_{23}
\end{equation}
\begin{equation}
b_{33}a'_{03}=a_{30}b_{12}^3+3a_{21}b_{12}^2b_{22}+3a_{12}b_{12}b_{22}^2+a_{03}b_{22}^3+3b_{12}b_{13}+3b_{22}b_{23}
\end{equation}
\end{lemma}

\noindent{\bf Proof:} This follows directly from $T_3$ and Equation 5.

We now come to the main result of this section, i.e. the standard forms of elliptical jets.

\begin{theorem}
An elliptical jet can be transformed by affine transformations to the following standard form \begin{equation}z(x,y)=\frac{1}{2}(x^2+y^2)+\frac{1}{6}(x^3-3xy^2)+\cdots\end{equation}
\noindent Or \begin{equation}z(x,y)=\frac{1}{2}(x^2+y^2)+\ \ \ \ \ \ \ \ 0 \ \ \ \ \ \ \ \ \ +\cdots\end{equation}\end{theorem}

We call the two cases the regular case and the degenerate case respectively.

\noindent {\bf Proof: }By Corollary 3.1, we can assume the elliptical jet is of the form $$z(x,y)=\frac{1}{2}(x^2+y^2)+\frac{1}{6}(a_{30}x^3+3a_{21}x^2 y+3a_{12}xy^2+a_{03}y^3)+\cdots$$
Let $T_4\in G_2$ be the transformation given by $$x=\sqrt{b_{33}}x'\cos t +\sqrt{b_{33}}y'\sin t+b_{13}z', x=-\sqrt{b_{33}}x'\sin t +\sqrt{b_{33}}y'\cos t+b_{23}z', z=b_{33}z'.$$ \noindent
By setting $a'_{12}=-a'_{30},a'_{21}=a'_{03}$ and using Equations 7 and 9, Equations 8 and 10, we obtain two linear equations with unknowns $b_{13}$ and $b_{23}$. Thus $b_{13},b_{23}$ can be determined. This means that we can transform a jet into a form $$z'(x',y')=\frac{1}{2}(x'^2+y'^2)+\frac{1}{6}(a'_{30}(x'^3-3x'y'^2)+a'_{03}(3x'^2 y'+y'^3))+\cdots.$$
By substituting $b_{13}, b_{23}$ into the Equations 7 and 10, we obtain
\begin{equation}a'_{30}=\sqrt{b_{33}}(\frac{1}{4}a_{30}\cos3t-\frac{3}{4}a_{21}\sin3t-\frac{3}{4}a_{12}\cos3t+\frac{1}{4}a_{03}\sin3t)\end{equation}
\begin{equation}a'_{03}=\sqrt{b_{33}}(\frac{1}{2}a_{30}\sin3t+\frac{3}{2}a_{21}\cos3t-\frac{3}{2}a_{12}\sin3t-\frac{1}{2}a_{03}\cos3t)\end{equation}
If we already have $a_{12}=-a_{30}$ and $a_{21}=a_{03}$, then Equations 13 and 14 become to
 $$a'_{30}=\sqrt{b_{33}}(a_{30}\cos3t-\frac{a_{03}}{2}\sin3t),\frac{a'_{03}}{2}=\sqrt{b_{33}}(a_{30}\sin3t+\frac{a_{03}}{2}\cos3t).$$
 It is obvious that besides the case both $a_{30}$ and $a_{03}$ are zero we can choose $t$ and $b_{33}>0$ suitably such that $a'_{30}=1$ and $a'_{03}=0$. This proves the theorem.

\begin{coro}
Let $G_3$ be the subgroup of $G_2$ which fixes the jets with $a_{00}=a_{10}=a_{01}=a_{11}=a_{21}=a_{03}=0,a_{20}=a_{02}=a_{30}=-a_{12}=1$, then $G_3$ is the symmetry group $D_3$ of equilateral triangle with vertices $(1,0),(-\frac{1}{2},\frac{\sqrt{3}}{2})$ and $(-\frac{1}{2},-\frac{\sqrt{3}}{2})$.
\end{coro}

\noindent{\bf Proof: }For affine transformation $T\in G_3$, by Equations 7, 8, 9 and 10, we have
%if $A=\left(\begin{array}{cc}\cos(t) & -\sin(t) \\\pm \sin(t) & \pm \cos(t) \\\end{array}\right)$ fixes the standard form \ref{standard}, then we have
$$b_{33}=b_{11}^3-3b_{11}b_{21}^2+3b_{11}b_{13}+3b_{21}b_{23}$$
$$0=b_{11}^2b_{12}-2b_{11}b_{21}b_{22}-b_{12}b_{21}^2+b_{12}b_{13}+b_{22}b_{23}$$
$$-b_{33}=b_{11}b_{12}^2-b_{11}b_{22}^2-2b_{12}b_{21}b_{22}+b_{11}b_{13}+b_{21}b_{23}$$
$$0=b_{12}^3-3b_{12}b_{22}^2+3b_{12}b_{13}+3b_{22}b_{23}$$
\noindent By canceling $b_{13}$ and $b_{23}$, we get

\begin{equation}4b_{33}=b_{11}^3-3b_{11}b_{12}^2-3b_{11}b_{21}^2+3b_{11}b_{22}^2+6b_{12}b_{21}b_{22}\end{equation}
\begin{equation}0=-3b_{11}^2b_{12}+6b_{11}b_{21}b_{22}+b_{12}^3+3b_{12}b_{21}^2-3b_{12}b_{22}^2\end{equation}

\noindent Since the transformation $T$ is in $G_2$, we must have $\left(
      \begin{array}{cc}
        b_{11} & b_{12} \\
        b_{21} & b_{22} \\
      \end{array}
    \right)=\pm\sqrt{b_{33}}\left(
                                                 \begin{array}{cc}
                                                   \cos t & -\sin t \\
                                                   \pm \sin t & \pm \cos t \\
                                                 \end{array}
                                               \right)$.
    By substituting $b_{11}, b_{12}, b_{21}, b_{22}$ into Equations 15 and 16, we have $b_{33}=1,\cos(3t)=1$. And we get the six solutions
$$\left(
      \begin{array}{cc}
        b_{11} & b_{12} \\
        b_{21} & b_{22} \\
      \end{array}
    \right)=\left(
    \begin{array}{cc}   1 & 0 \\   0 & \pm 1 \\ \end{array}
\right),\left(    \begin{array}{cc}
\cos\frac{2}{3}\pi & -\sin\frac{2}{3}\pi \\
      \pm \sin\frac{2}{3}\pi & \pm \cos\frac{2}{3}\pi \\
                                                 \end{array}
                                               \right),\left(
                                                 \begin{array}{cc}
                                                   \cos\frac{4}{3}\pi & -\sin\frac{4}{3}\pi \\
                                                   \pm \sin\frac{4}{3}\pi & \pm \cos\frac{4}{3}\pi \\
                                                 \end{array}
                                               \right).$$

\noindent The generators of $D_3$ are the $\frac{2\pi}{3}$ rotation $\sigma$ and the involution $\tau$ changing the sign of coordinate $y$. Their action on coordinates $x,y,z$ are given by

\begin{equation}\sigma: x=-\frac{1}{2}x'-\frac{\sqrt{3}}{2}y',y=\frac{\sqrt{3}}{2}x'-\frac{1}{2}y',z=z'; \ \ \ \ \tau: x=x',y=-y',z=z'.\end{equation}

We summarize this section into an algorithm A1 to determine the standard form of a jet and the corresponding affine transformation.

{\bf The Algorithm A1:}

Step1.Start with a surface of local form $z=z(x,y)$, expand it at a point $(x_0,y_0)$ as Equation 1.

Step2. Using a translation $T_0$ in $Aff(3)$ to transform the jet into the form of Equation 2.

Step3. Using affine transformation $T_1$ preserving $(0,0,0)$ to transform the jet into the form of Equation 4.

Step4. Using affine transformation $T_2$ in $G_1\subset Aff(3)$ (see Lemma3.2) to transform a jet of elliptical surface into the form of Equation 6.

Step5. Using affine transformation $T_3\in G_2\subset Aff(3)$(see Lemma3.4) to transform a jet of elliptical surface into the form of Equation 11 or 12.

Now we have the output of Algorithm A1 of the standard form of a regular elliptical jet and an affine transformation $T=T_3T_2T_1T_0$.

\section{Equivariant moving frame and its moving equations}
In this section, we only discuss the regular elliptical surfaces.
%Using jet space as a tool, we can transform the study of local affine congruence invariants of surfaces $S$ to the invariants of the jet of surfaces $S$. We have the following definition.

\subsection{Differential invariants and invariant forms}

\begin{define}
A smooth function $f:J^{2,r}(A^3)\to \mathbb{R}$ which is invariant under the action of $Aff(3)$ is called a differential invariant of surfaces in $A^3$ of order not great than $r$.
\end{define}

%For a surfaces $S$ in $A^3$. We can regard $S$ as the surface $j(S)$ in jet space $J^{2,r}(A^3)$.
Up to affine congruence the jet of a regular elliptical surface $S$ at a point $P$ can be written as a form %$z=\sum\limits_{n=0}^\infty \frac{1}{n!}(\sum\limits_{i,j\geq 0}\choose{n}{i} a_{ij}(x-x_0)^i(y-y_0))$, its standard form
$$\D{z=\frac{1}{2}(x^2+y^2)+\frac{1}{6}(x^3-3xy^2)+\sum\limits_{n=4}^\infty \frac{1}{n!}(\sum\limits_{i+j=n} I_{ij}{{n}\choose{i}}x^i y^j)}.$$ It is obvious that the coefficients $I_{ij},i+j\geq 4$ are all affine differential invariants on surface $S$. Recall that in Definition 2.2, we defined the invariantization $\iota(F)$ of a function $F$. In fact, $I_{ij}$ is the invariantization of $a_{ij}$, i.e. $\iota(a_{ij})=I_{ij}$. The 1-forms $dx,dy$ can also be invariantized into affine invariant differential form $w_1,w_2$. The corresponding dual invariant differential operators $\mathrm{D}^1,\mathrm{D}^2$ are defined by $df=\mathrm{D}^1f\  w_1+\mathrm{D}^2f\  w_2$. For details see \cite{Fels_Olver_04}.

 By Theorem 3.1, the jet of $S$ at a point $(x_0,y_0,z_0)$ can be transformed to the standard form by an affine transformation $T$. Let the transformation $T$ be of the form $$\left(
\begin{array}{c}
 x \\
  y \\
  z\\
  1 \\
   \end{array}
   \right)=\left(
             \begin{array}{cccc}
               b_{11} & b_{12} & b_{13} & x_0\\
               b_{21} & b_{22} & b_{23} & y_0\\
               b_{31} & b_{32} & b_{33} & z_0\\
               0&0&0&1\\
             \end{array}
           \right)\left(
\begin{array}{c}
 x' \\
  y' \\
  z'\\
   1\\
   \end{array}
   \right)$$
\noindent We set $e_1=\left(
                   \begin{array}{c}
                     b_{11} \\
                     b_{21} \\
                     b_{31} \\
                   \end{array}
                 \right)
,e_2=\left(
                   \begin{array}{c}
                     b_{12} \\
                     b_{22} \\
                     b_{32} \\
                   \end{array}
                 \right),e_3=\left(
                   \begin{array}{c}
                     b_{13} \\
                     b_{23} \\
                     b_{33} \\
                   \end{array}
                 \right),r=\left(
                   \begin{array}{c}
                     x_0 \\
                     y_0 \\
                     z_0 \\
                   \end{array}
                 \right)$. $(e_1,e_2,e_3,r)$ gives a smooth affine moving frame on surface $S$. The moving equation of the moving frame has the form
\begin{equation}d(e_1,e_2,e_3,r)=(e_1,e_2,e_3,r)\Omega \end{equation}
\noindent Where $\Omega$ is a matrix of 1-forms.

Since we have chosen $a_{00}=a_{10}=a_{01}=0$, which means that the fundamental vector fields $e_1,e_2$ dual to $w_1,w_2$ are tangent to $S$ and $dr=w_1e_1+w_2e_2$.
So $\Omega$ must be of the form
\begin{equation}
\Omega=\left(
                       \begin{array}{cccc}
                         w_{11} & w_{12} &  w_{13} &w_1 \\
                         w_{21} & w_{22} &  w_{23} &w_2 \\
                         w_{31} & w_{32} &  w_{33} &0   \\
                         0&0&0&0\\
                       \end{array}
                     \right)\end{equation}
\noindent The Maurer-Cartan invariants $R_{ij}^k,1\leq i,j\leq 3,1\leq k\leq 2$ are defined by
\begin{equation}w_{ij}=\sum\limits_{k=1}^2 R_{ij}^k w_k.\end{equation}

For differential manifold $M$, let $Diff(M)$ be its diffeomorphism group and $\mathcal{X}(M)$ be the Lie algebra of smooth vector fields. The action of $Aff(3)$ on $A^3$ can be written as a homomorphism $Aff(3)\stackrel{\Phi}\to Diff(A^3)$. ${\Phi}$ induces the Lie algebra homomorphism $aff(3)\stackrel{\phi} \rightarrow \mathcal{X}(A^3)$. The Lie algebra $aff(3)$ has a natural basis $X_i, X_{ij},1\leq i,j\leq 3$. Under the homomorphism $\phi$ the basis corresponds to vector fields

$$\D{X_1=\frac{\partial}{\partial x}, X_2=\frac{\partial}{\partial y},X_3=\frac{\partial}{\partial z}},$$
$$\D{X_{11}=x\frac{\partial}{\partial x}, X_{12}=y\frac{\partial}{\partial x}, X_{13}=z\frac{\partial}{\partial x}}$$
$$\D{X_{21}=x\frac{\partial}{\partial y}, X_{22}=y\frac{\partial}{\partial y}, X_{23}=z\frac{\partial}{\partial y}},$$
$$\D{X_{31}=x\frac{\partial}{\partial z}, X_{32}=y\frac{\partial}{\partial z}, X_{33}=z\frac{\partial}{\partial z}}.$$

Similarly the induced action of $Aff(3)$ on $J^{2,3}(A^3)$ can be written as a homomorphism $Aff(3)\to Diff(J^{2,3}(A^3))$. It induces the Lie algebra homomorphism $aff(3)\stackrel{\psi} \rightarrow \mathcal{X}(J^{2,3}(A^3))$. Under the homomorphism $\psi$, the vector fields corresponding to the basis $X_i, X_{ij},1\leq i,j\leq 3$ are the prolongation of the above vector fields. These vector fields can be calculated by the Theorem 4.16 in Olver\cite{Olver_95}.

Computation shows that on $J^{2,3}(A^3)$ these vector fields are given by

$$\D{X_1=\frac{\partial}{\partial x}, X_2=\frac{\partial}{\partial y},X_3=\frac{\partial}{\partial a_{00}}},$$
$$\D{X_{11}=x\frac{\partial}{\partial x}-a_{10}\frac{\partial}{\partial a_{10}}-2a_{20}\frac{\partial}{\partial a_{20}}-a_{11}\frac{\partial}{\partial a_{11}}-3a_{30}\frac{\partial}{\partial a_{30}}-2a_{21}\frac{\partial}{\partial a_{21}}-a_{12}\frac{\partial}{\partial a_{12}}},$$
$$\D{X_{12}=y\frac{\partial}{\partial x}-a_{10}\frac{\partial}{\partial a_{01}}-a_{20}\frac{\partial}{\partial a_{11}}-2a_{11}\frac{\partial}{\partial a_{02}}-a_{30}\frac{\partial}{\partial a_{21}}-2a_{21}\frac{\partial}{\partial a_{12}}-3a_{12}\frac{\partial}{\partial a_{03}}},$$
$$\D{X_{13}=a_{00}\frac{\partial}{\partial x}-a^2_{10}\frac{\partial}{\partial a_{10}}-a_{10} a_{01}\frac{\partial}{\partial a_{01}}-3 a_{10} a_{20}\frac{\partial}{\partial a_{20}}-(a_{01} a_{20} + 2 a_{10} a_{11})\frac{\partial}{\partial a_{11}}-(2 a_{01} a_{11} +a_{10} a_{02})\frac{\partial}{\partial a_{02}}}$$
$$\D{-(4 a_{10} a_{30} + 3 a^2_{20})\frac{\partial}{\partial a_{30}}-(a_{01} a_{30}+ 3 a_{10} a_{21} + 3 a_{20}a_{11})\frac{\partial}{\partial a_{21}}-(2 a_{01} a_{21} +a_{20} a_{02}+ 2 a_{10} a_{12} + 2 a^2_{11})\frac{\partial}{\partial a_{12}}}$$
$$\D{+(3 a_{01} a_{12} + 3 a_{11} a_{02} +a_{10} a_{03})\frac{\partial}{\partial a_{03}}},$$
$$\D{X_{21}=x\frac{\partial}{\partial y}-a_{01}\frac{\partial}{\partial a_{10}}-2a_{11}\frac{\partial}{\partial a_{20}} -a_{02}\frac{\partial}{\partial a_{11}}-3a_{21}\frac{\partial}{\partial a_{30}}-2a_{12}\frac{\partial}{\partial a_{21}}-a_{03}\frac{\partial}{\partial a_{12}}},$$
$$\D{X_{22}=y\frac{\partial}{\partial y}-a_{01}\frac{\partial}{\partial a_{01}}-a_{11}\frac{\partial}{\partial a_{11}}-2 a_{02}\frac{\partial}{\partial a_{02}}-a_{21}\frac{\partial}{\partial a_{21}}-2 a_{12}\frac{\partial}{\partial a_{12}}-3 a_{03}\frac{\partial}{\partial a_{03}}},$$
$$\D{X_{23}=a_{00}\frac{\partial}{\partial y}-a_{01} a_{10}\frac{\partial}{\partial a_{10}}-a^2_{01}\frac{\partial}{\partial a_{01}}-(a_{01} a_{20} + 2 a_{10} a_{11})\frac{\partial}{\partial a_{20}}-(2 a_{01} a_{11} +a_{10} a_{02})\frac{\partial}{\partial a_{11}}-3 a_{01} a_{02}\frac{\partial}{\partial a_{02}}}$$
$$\D{-(a_{01} a_{30} + 3 a_{10} a_{21} + 3 a_{20} a_{11} ) \frac{\partial}{\partial a_{30}}-(2 a_{01} a_{21} +a_{20} a_{02} + 2 a_{10} a_{12}
+ 2 a^2_{11})\frac{\partial}{\partial a_{21}}}$$
$$\D{-(3 a_{01} a_{12} + 3 a_{11} a_{02} +a_{10} a_{03})\frac{\partial}{\partial a_{12}}-(4 a_{01} a_{03} + 3 a^2_{02}) \frac{\partial}{\partial a_{03}}},$$
$$\D{X_{31}=x\frac{\partial}{\partial a_{00}}+\frac{\partial}{\partial a_{10}}},$$
$$\D{X_{32}=y\frac{\partial}{\partial a_{00}}+\frac{\partial}{\partial a_{01}}},$$
$$\D{X_{33}=a_{00}\frac{\partial}{\partial a_{00}}+a_{10}\frac{\partial}{\partial a_{10}}+a_{01}\frac{\partial}{\partial a_{01}}+a_{20}\frac{\partial}{\partial a_{20}}+a_{11}\frac{\partial}{\partial a_{11}}+a_{02}\frac{\partial}{\partial a_{02}}}$$
$$\D{+a_{30}\frac{\partial}{\partial a_{30}}+a_{21}\frac{\partial}{\partial a_{21}}+a_{12}\frac{\partial}{\partial a_{12}}+a_{03}\frac{\partial}{\partial a_{03}}}.\ \ \ \ \ \ \ \ \ \ \ \  \ (*)$$

The Maurer-Cartan invariants can be computed by the following theorem of Fels and Olver\cite{Fels_Olver_98}\cite{Fels_Olver_99}.
\begin{theorem}
Let $F$ be a differential function on jet spaces $J^{2,r}(A^3)$ and $\iota(F)$ be its invariantization. Then
$$\mathrm{D}^k(\iota(F))=\iota(D^k(F))-\sum\limits_{i,j=1}^3 \iota(X_{ij}F)R_{ij}^k.$$
\noindent Where $R_{ij}^k,1\leq i,j\leq 3,1\leq k\leq 2$ are the Maurer-Cartan differential invariants.
\end{theorem}
In this theorem, $D^k$ represents the partial differential operator. For example $D^1(a_{21})=a_{31}$ and $D^2(a_{21})=a_{22}$.
%Procedure 2: Compute the Maurer-Cartan invariants from the fixed cross section.
%$X_{12}=y\frac{\partial}{\partial x}+ \frac{\partial}{\partial y}-\frac{\partial}{\partial a_{00}}+\frac{\partial}{\partial a_{10}}+\frac{\partial}{\partial a_{01}}+\frac{\partial}{\partial a_{20}}+\frac{\partial}{\partial a_{11}}+\frac{\partial}{\partial a_{02}}+\frac{\partial}{\partial a_{30}}+\frac{\partial}{\partial a_{21}}+\frac{\partial}{\partial a_{12}}\frac{\partial}{\partial a_{03}}$

\begin{lemma}The Maurer-Cartan invariants are
$$\D{\left(
  \begin{array}{ccc}
    R_{11}^1 & R_{12}^1 & R_{13}^1 \\
    R_{21}^1 & R_{22}^1 & R_{23}^1 \\
    R_{31}^1 & R_{32}^1 & R_{33}^1 \\
  \end{array}
\right)=\left(
          \begin{array}{ccc}
            -\frac{I_{40}}{4}+\frac{3I_{22}}{4}-\frac{1}{2} & -\frac{I_{31}}{4}+\frac{I_{13}}{12} & -\frac{I_{40}}{4}-\frac{I_{22}}{4}+\frac{1}{2} \\
            \frac{I_{31}}{4}-\frac{I_{13}}{12} & -\frac{I_{40}}{4}+\frac{3I_{22}}{4}+\frac{1}{2} & -\frac{I_{31}}{4}-\frac{I_{13}}{4} \\
            1 & 0 & -\frac{I_{40}}{2}+\frac{3I_{22}}{2} \\
          \end{array} \right)}$$
$$\D{\left(
  \begin{array}{ccc}
    R_{11}^2 & R_{12}^2 & R_{13}^2 \\
    R_{21}^2 & R_{22}^2 & R_{23}^2 \\
    R_{31}^2 & R_{32}^2 & R_{33}^2 \\
  \end{array}
\right)=\left(
          \begin{array}{ccc}
            -\frac{I_{31}}{4}+\frac{3I_{13}}{4} & -\frac{I_{22}}{4}+\frac{I_{04}}{12}+\frac{1}{2} & -\frac{I_{31}}{4}-\frac{I_{13}}{4} \\
            \frac{I_{22}}{4}-\frac{I_{04}}{12}+\frac{1}{2} & -\frac{I_{31}}{4}+\frac{3I_{13}}{4} & -\frac{I_{22}}{4}-\frac{I_{04}}{4}+\frac{1}{2} \\
            0 & 1 & -\frac{I_{31}}{2}+\frac{3I_{13}}{2} \\
          \end{array}
        \right)}$$
\end{lemma}

\noindent {\bf Proof: }By the expressions of vector fields $X_{ij}$ and Theorem 4.1, under the invariantization condition $a_{10}=a_{01}=a_{11}=0, a_{20}=a_{02}=1$, we have

$R_{31}^1-\iota(a_{20})=0,\ \ \ \ \ \ R_{32}^2-\iota(a_{11})=0$

$R_{32}^1-\iota(a_{11})=0,\ \ \ \ \ \ R_{32}^2-\iota(a_{02})=0$

$-2R_{11}^1+R_{33}^1-\iota(a_{30})=0,\ \ \ \ \ \ -2R_{11}^2+R_{33}^2-\iota(a_{21})=0$

$-R_{12}^1-R_{21}^1-\iota(a_{21})=0,\ \ \ \ \ \ -R_{12}^2-R_{21}^2-\iota(a_{12})=0$

$-2R_{22}^1+R_{33}^1-\iota(a_{12})=0,\ \ \ \ \ \ -2R_{22}^2+R_{33}^2-\iota(a_{03})=0$

Solving these equations we have

$\D{R_{31}^1=1, R_{32}^1=0, R_{11}^1=\frac{R_{33}^1-1}{2}, R_{22}^1=\frac{R_{33}^1+1}{2},R_{21}^1=-R_{12}^1;}$

$\D{R_{31}^2=0, R_{32}^2=1, R_{11}^2=\frac{R_{33}^2}{2}, R_{22}^2=\frac{R_{33}^2}{2},R_{21}^2=-R_{12}^2+1. }$

The conditions $a_{30}=-a_{12}=1$, and $a_{21}=a_{03}=0$ give the equations

$-3R_{11}^1-3R_{13}^1+R_{33}^1-\iota(a_{40})=0$

$-R_{12}^1+2R_{21}^1-R_{23}^1-\iota(a_{31})=0$

$R_{11}^1-R_{13}^1+2R_{22}^1-R_{33}^1-\iota(a_{22})=0$

$3R_{12}^1-3R_{23}^1-\iota(a_{13})=0$

and

$-3R_{11}^2-3R_{13}^2+R_{33}^2-\iota(a_{31})=0$

$-R_{12}^2+2R_{21}^2-R_{23}^2-\iota(a_{22})=0$

$R_{11}^2-R_{13}^2+2R_{22}^2-R_{33}^2-\iota(a_{13})=0$

$3R_{12}^2-3R_{23}^2-\iota(a_{04})=0$

Solving these equations , we prove the lemma.

%$R_{23}^1=-\frac{I_{31}}{4}-\frac{I_{13}}{4},  R_{12}^1=-\frac{I_{31}}{4}+\frac{I_{13}}{12}, R_{21}^1=\frac{I_{31}}{4}-\frac{I_{13}}{12}, R_{33}^1=-\frac{I_{40}}{2}+\frac{3I_{22}}{2}$.

%$R_{13}^1=-\frac{I_{40}}{4}-\frac{I_{22}}{4}+\frac{1}{2},R_{11}^1=-\frac{I_{40}}{4}+\frac{3I_{22}}{4}-\frac{1}{2},R_{22}^1=-\frac{I_{40}}{4}+\frac{3I_{22}}{4}+\frac{1}{2}$

%$R_{23}^1=-\frac{I_{13}+I_{31}}{4},$

\begin{theorem}
Let $S$ be a regular elliptical surface, then the 1-forms $w_1,w_2$ and $I_{40}, I_{31},I_{22}$ are complete affine differential invariants.
\end{theorem}

\noindent {\bf Proof: }By Cartan\cite{Cartan_37}, the moving equations of the moving frame determine the surface up to congruence. Hence by Equations 18, 19 and Lemma 4.1, $w_1,w_2$ and $I_{40}, I_{31}, I_{22}, I_{13}, I_{04}$ are complete differential invariants. Define invariants $Y_1,Y_2$ by $dw_1=Y_1 w_1\wedge w_2, dw_2=Y_2 w_1\wedge w_2$. Computation shows $\D{Y_1=R_{11}^2-R_{12}^1=\frac{2I_{13}}{3}}$ and $\D{Y_2=R_{21}^2-R_{22}^1=\frac{I_{40}}{4}-\frac{I_{22}}{2}-\frac{I_{04}}{12}}$.
Since $Y_1,Y_2$ are invariants determined by $w_1,w_2$, and $I_{13}, I_{04}$ can be computed from $Y_1,Y_2$ and $I_{40},I_{31}, I_{22}$, this proves the theorem.

\begin{define}
We call the invariants $I_{40},I_{31},I_{22},I_{13},I_{04}$ the fundamental curvature invariants of elliptical surfaces.
\end{define}

\begin{define}

The first fundamental forms of a regular elliptical surface are defined to be the 1-forms $\Phi^1_1=w_1$ and $\Phi^2_1=w_2$. The second, third and fourth fundamental forms of a regular elliptical surface are $\Phi_2=w_1^2+w_2^2, \Phi_3=w_1^3-3w_1w_2^2, \Phi_4=I_{40}w_1^4+4I_{31}w_1^3 w_2+ 6I_{22}w_1^2 w_2^2+4I_{13}w_1 w_2^3+I_{04}w_2^4$.
\end{define}

\begin{theorem}
Let $S$ be a regular elliptical surface, then the second, third and fourth fundamental forms $\Phi_2,\Phi_3,\Phi_4$ are complete invariants.
\end{theorem}

\noindent {\bf Proof: }We note that if $w_1,w_2$ is a solution of equations $w_1^2+w_2^2=\Phi_2, w_1^3-3w_1w_2^2=\Phi_3$, then other solutions are given by
the action of symmetry group $D_3$ in Corollary 3.2. Hence we can solve $w_1,w_2$ from $\Phi_2,\Phi_3$. By expanding $\Phi_4$ with respect to $w_1,w_2$, we can determine $I_{40}, I_{31}, I_{22}, I_{13}, I_{04}$. So $\Phi_2,\Phi_3,\Phi_4$ determine $w_1,w_2$ and $I_{40}, I_{31}, I_{22}, I_{13}, I_{04}$ up to the action of $D_3$. This proves the theorem.

\subsection{Compatible conditions}
By Theorem 4.2, the 1-forms $w_1,w_2$ and $I_{40}, I_{31},I_{22},I_{13}, I_{04}$ are complete differential invariants. But not any set of $w_1,w_2$ and $I_{40}, I_{31},I_{22},I_{13}, I_{04}$ can be realized as the invariants of some surfaces. They must satisfy some compatible conditions which are given by the Cartan's structure equation $d\Omega+\Omega\wedge \Omega=0$. The components of this equation are

\noindent 1. $dw_{ij}=-w_{ik}\wedge w_{kj}$, $1\leq i,j\leq 3$. %\sum\limits_{k=1}^{3}

\noindent 2. $dw_1=-w_{11}\wedge w_1-w_{12}\wedge w_2, dw_2=-w_{21}\wedge w_1-w_{22}\wedge w_2.$

\noindent For $dw_{ij}=-w_{ik}\wedge w_{kj}$, the left side of this equation is

\noindent $dw_{ij}=d(R_{ij}^k w_k)=d(R_{ij}^k)\wedge w_k+R_{ij}^k dw_k$

$=\mathrm{D}^1R_{ij}^k w_1\wedge w_{k}+\mathrm{D}^2 R_{ij}^k w_2\wedge w_k+R_{ij}^k Y_kw_1\wedge w_2$.

 $=(\mathrm{D}^1R_{ij}^2-\mathrm{D}^2 R_{ij}^1+R_{ij}^1 Y_1+R_{ij}^2 Y_2)w_1\wedge w_2$.

\noindent The right side of the equation is

\noindent $-w_{ik}\wedge w_{kj}=-R_{ik}^l w_l \wedge R_{kj}^m w_m= (R_{ik}^2 R_{kj}^1-R_{ik}^1 R_{kj}^2)w_1\wedge w_2$.

\noindent Hence the structure equations give
\begin{equation} \mathrm{D}^1R_{ij}^2-\mathrm{D}^2 R_{ij}^1+R_{ij}^1 Y_1+R_{ij}^2 Y_2=R_{ik}^2 R_{kj}^1-R_{ik}^1 R_{kj}^2.\end{equation}

%\noindent Since $w_{31}=w_1, w_{32}=w_2$, we have $R_{31}^1=1=R_{32}^2=1, R_{31}^2=1=R_{32}^1=0$.

%The structure equations is $dw_1=Y_1w_1\wedge w_2, dw_2=Y_2w_1\wedge w_2, 0=$

\begin{theorem}
The compatible conditions are given by

\noindent $ 36\mathrm{D}^2I_{31}-36\mathrm{D}^1I_{22}-12\mathrm{D}^2I_{13}+12\mathrm{D}^1I_{04}-9I_{40}I_{22}+3I_{40}I_{04}-24I_{31}I_{13}+18I_{22}^2-3I_{22}I_{04}+8I_{13}^2-I_{04}^2-18I_{40}-36I_{22}-18I_{04}=0.$

\noindent $ 12\mathrm{D}^2I_{40}-12\mathrm{D}^1I_{31}+12\mathrm{D}^2I_{22}-12\mathrm{D}^1I_{13}-3I_{40}I_{31}-23I_{40}I_{13}+18I_{31}I_{22}+5I_{31}I_{04}-6I_{22}I_{13}+I_{13}I_{04}+48I_{13}=0.$

\noindent $ 12\mathrm{D}^2I_{31}-12\mathrm{D}^1I_{22}+12\mathrm{D}^2I_{13}-12\mathrm{D}^1I_{04}-3I_{40}I_{22}-7I_{40}I_{04}-16I_{31}I_{13}+18I_{22}^2+5I_{22}I_{04}-16I_{13}^2+I_{04}^2+18I_{40}-36I_{22}-6I_{04}=0.$

\noindent $12\mathrm{D}^2I_{40}+12\mathrm{D}^1 I_{31}-36\mathrm{D}^2I_{22}+36\mathrm{D}^1I_{13}-3I_{40}I_{31}+I_{40}I_{13}+6I_{31}I_{22}+I_{31}I_{04}+6I_{22}I_{13}-3I_{13}I_{04}=0.$

\noindent $dw_{1}=Y_1 w_1\wedge w_2.$

\noindent $dw_{2}=Y_2 w_1\wedge w_2.$
\end{theorem}

\noindent {\bf Proof}: This is computed directly from Equation 21 and Lemma 4.1. The first compatible condition corresponds to $i=1, j=2$ in Equation 21, the second corresponds to $i=1, j=3$, the third corresponds to $i=2, j=3$,
the fourth corresponds to $i=3, j=3$, the fifth corresponds to $i=3, j=1$ and the sixth corresponds to $i=3, j=2$.

\begin{coro}
Let $S$ be a regular elliptical surface with constant curvatures $I_{40},I_{31},I_{22},I_{13},I_{04}$, then these invariants satisfy
\begin{equation}( -23I_{40}-6I_{22}+I_{04}+48) I_{13}+ ( -3I_{40}+18I_{22}+5I_{04}) I_{31}=0\end{equation}
\begin{equation}( I_{40}+6I_{22}-3I_{04})I_{13}+ (-3I_{40}+6I_{22}+I_{04}) I_{31}=0\end{equation}
\begin{equation}-9I_{40}I_{22}+3I_{40}I_{04}-24I_{31}I_{13}+18I_{22}^2-3I_{22}I_{04}+8I_{13}^2-I_{04}^2-18I_{40}-36I_{22}-18I_{04}=0\end{equation}
\begin{equation}-3I_{40}I_{22}-7I_{40}I_{04}-16I_{31}I_{13}+18I_{22}^2+15I_{22}I_{04}-16I_{13}^2+I_{04}^2+18I_{40}-36I_{22}-6I_{04}=0\end{equation}
\end{coro}

\noindent {\bf Proof}: The Equations 22, 23, 24 and 25 are derived directly from the second, the fourth, the first and the third equations of Theorem 4.4 respectively.

Now we can summarize our second algorithm A2 which computes the moving equations and compatible conditions.

{\bf The Algorithm A2:}

Step1. Start from the output of the Algorithm A1. i.e. a standard form of a regular elliptical surface $S$ at a point $P$ and a corresponding affine transformation $T$.
$T$ gives the moving frame.

%Step2. Compute the moving equation (18), here the coefficient matrix $\Omega$ has the form of (19).
%The Maurer-Cartan invariants $R^k_{ij}$ are defined by (20).

Step2. Compute the induced infinitesimal action of $Aff(3)$ on jet space $J^{2,3}(A^3)$ by (*). We compute it with the Maple package ¡®JetCalculus¡¯

Step3. Compute the Maurer-Cartan invariants $R^k_{ij}$ by Theorem 4.1.

Step4. Compute the compatible conditions from the Cartan's structure equation $d\Omega+\Omega\wedge \Omega=0$.

%The Maurer-Cartan invariants $R^k_ij$ is actually equivalent to the moving equation.
Once the Maurer-Cartan invariants are computed, the moving equations are determined. Combining the compatible conditions, we have the main ingredient of affine geometry of elliptical surfaces.

\section{Classification of elliptical surfaces with constant curvatures}

In this section we classify the elliptical surfaces with constant curvatures. We give an example at first.
\begin{example}

Considering the surface $S: xy^k z^k=1$, We compute its moving equation by algorithms A1 and A2.

For a point $P(x,y,z)\in S$ with $xyz\not=0$.
Using the affine transformation $T_1: x_1=x x_2, y_1=y y_2, z_1=z z_2$, we transform the jet of $S$ at $P$ to the point $(1,1,1)$.

By the transformation $T_2: x_2=x_3+1,y_2=y_3+1,z_2=z_3+1$, the jet is transformed to a jet at $(0,0,0)$. Locally it has the form

$z_3=\frac{1}{2}(k(k+1)x_3^2+k^2x_3y_3+k(k+1)y_3^2)$

$+\frac{1}{6}(-k(k+1)(k+2)x_3^3+3k^2(k+1)x_3^2y_3+3k^2(k+1)x_3y_3^2-k(k+1)(k+2)y_3^3)$

$+\frac{1}{24}(k(k+1)(k+2)(k+3)x_3^4+4k^2(k+1)(k+2)x_3^3y_3+6k^2(k+1)^2x_3^2y_3^2+4k^2(k+1)(k+2)x_3y_3^3$

$+k(k+1)(k+2)(k+3)y_3^4)+\cdots$.

The quadratic part is positively definite. So the surface is elliptical.

By the transformation $T_3: x_3=\frac{\sqrt{2}}{2\sqrt{2k^2+k}}x_4-\frac{\sqrt{2}}{2\sqrt{k}}y_4, y_3=\frac{\sqrt{2}}{2\sqrt{2k^2+k}}x_4+\frac{\sqrt{2}}{2\sqrt{k}}y_4, z_3=z_4$, the jet is transformed to

$z_4=\frac{1}{2}(x_4^2+y_4^2)+\frac{1}{6}\frac{\sqrt{2}(k+1)}{\sqrt{k(2k+1)}})(x_4^3+3x_4y_4^2)+\frac{1}{24}(\frac{(2k+3)(k+1)}{k(2k+1)}x_4^4+6\frac{(2k+3)(k+1)}{k(2k+1)}x_4^2 y_4^2+\frac{(6k+3)(k+1)}{k(2k+1)}y_4^4+\cdots$.

By the transformation $T_4: x_4=\frac{\sqrt{2 k(2k+1)}}{k+1}x_5+\frac{\sqrt{2 k(2k+1)}}{k+1}z_5,y_4=-\frac{\sqrt{2 k(2k+1)}}{k+1}y_5,z_4=\frac{2k(2k+1)}{{(k+1)}^2} z_5$, the jet can be transformed to a jet of local form

$$z_5=\frac{1}{2}(x_5^2+y_5^2)+\frac{1}{6}(x_5^3-3x_5y_5^2)+\frac{1}{24}(-\frac{(k-1)}{k+1}x_5^4-6 \frac{(k-1)}{k+1}x_5^2y_5^2+\frac{3(k-1)}{k+1}y_5^4)+\cdots$$

The composition $T=T_4T_3T_2T_1$ is given by

\begin{center}$x_1=\frac{x}{k+1}x_5+\frac{\sqrt{2k+1}y}{k+1}y_5+\frac{z}{k+1}z_5+x$

$y_1=\frac{x}{k+1}x_5-\frac{\sqrt{2k+1}y}{k+1}y_5+\frac{z}{k+1}z_5+y$

$z_1==-\frac{2kx}{k+1}x_5+\frac{2k^2z}{{(k+1)}^2} z_5+z$
\end{center}
The moving frames on $S$ are given by
$$r=(x,y,z)^t, e_1=\frac{1}{(k+1)}(x,y,-2kz)^t,e_2=\frac{\sqrt{2k+1}}{{k+1}}(x,-y,0)^t,e_3=\frac{1}{k+1}(x,y,\frac{2k^2}{k+1}z)^t.$$ Therefore if we let $A=\left( \begin {array}{cccc} {\frac {x}{k+1}}&{\frac {\sqrt {2\,k+1}x}{
k+1}}&{\frac {x}{k+1}} &x \\ \noalign{\medskip}{\frac {y}{k+1}}&-{\frac {
\sqrt {2\,k+1}y}{k+1}}&{\frac {y}{k+1}}& y\\ \noalign{\medskip} {
\frac {-2kz}{k+1}}&0& {\frac {2{k}^{2}z}{ \left( k+1 \right) ^{2}}}&z\\
0&0&0&1 \end {array} \right)$,
then

$\Omega=A^{-1}dA=\left( \begin {array}{cccc} -\frac{1}{2}\,{\frac { k\left( {\it y}\,dx+{\it x}
\,dy \right)}{xy}}&\frac{1}{2}\,{\frac {k \left( {\it y}\, dx-{\it x }\,dy
 \right) }{\sqrt {2\,k+1}xy}}&\frac{1}{2}\,{\frac { k\left( {\it y}\, dx+{\it
 x}\,dy \right)}{xy}}  & \frac{1}{2}\,{\frac {\left( k+1 \right) \left( {\it y}\, dx+{\it x}\,dy
 \right)   }{xy}} \\ \noalign{\medskip}\frac{1}{2}\,{\frac {{\it y}\, dx-{
\it x }\,dy}{\sqrt {2\,k+1}xy}}&\frac{1}{2}\,{\frac {{\it y}\,dx+{\it x}\,dy}{x
y}}&\frac{1}{2}\,{\frac {{\it y}\, dx-{\it x}\,dy}{\sqrt {2\,k+1}xy}} & \frac{1}{2}\,{\frac { \left( k+1 \right)
 \left( {\it y}\,dx-{\it x}\,dy \right) }{\sqrt {2\,k+1}xy}}
\\ \noalign{\medskip}\frac{1}{2}\,{\frac {\left( k+1 \right)  \left( {\it y}\,dx+{\it x}\,dy
 \right)  }{xy}}&\frac{1}{2}\,{\frac { \left( k+1 \right)
 \left( {\it y}\,dx-{\it x }\,dy \right) }{\sqrt {2\,k+1}xy}}&-\frac{1}{2}\,{\frac {(k-1) \left( {\it y}\,dx+{\it
 x}\,dy \right)}{xy}} & 0\\ 0&0&0&0
\end {array} \right)$

Hence we have $w_1=\frac{1}{2}\,{\frac {\left( k+1 \right) \left( {\it y}\, dx+{\it x}\,dy
 \right)   }{xy}}, w_2=\frac{1}{2}\,{\frac { \left( k+1 \right)
 \left( {\it y}\,dx-{\it x}\,dy \right) }{\sqrt {2\,k+1}xy}} $ and

 \begin{center}$\Omega=A^{-1}dA=\left(
                       \begin{array}{cccc}
                         -{\frac {k}{k+1}}w_1 &  {\frac {k}{k+1}}w_2 &  {\frac {k}{k+1}}w_1 &w_1\\
                         {\frac {1}{k+1}}w_2 &  {\frac {1}{k+1}}w_1 &  {\frac {1}{k+1}}w_2 & w_2\\
                         w_1 & w_2 & -\frac{k-1}{k+1}w_1 &0 \\ 0&0&0&0
                       \end{array}
                     \right)
$
\end{center}
As a result, we show that $S$ is a constant curvature surface with $I_{40}=I_{22}=-{\frac {k-1}{k+1}}, I_{04}={\frac {3(k-1)}{k+1}}, I_{31}=I_{13}=0$.
%define $n=e_3-R^1_{33}e_2-R^2_{33}e_2$, then we have a new moving frame $(e_1,e_2,n,r)$ and $w'_{33}=0$.
\end{example}

\begin{theorem}
%The solutions of Equation 19,20,21,22 must be of the form in the following table.
Let $S$ be a regular elliptical surface with constant curvatures, then  $I_{40},I_{31},I_{22},I_{13},I_{04}$ must be of the form in the following table.
\end{theorem}

$                                    \begin{array}{|c|c|c|c|c|c|}
                                     \hline    & I_{40} & I_{31} & I_{22} & I_{13} & I_{04}   \\
                                     \hline  1  & k & 0 & -1 & 0 & 3  \\
                                     \hline  2  & -3k^2+6k & 0 & -k^2-2k+2 & 0 & -3k^2-6k   \\
                                     \hline  3  & k & l(2k-3) & k & 0 & -3k   \\
                                     %\hline  4  & 9/4 & \mp 3\sqrt{3}/4 & -1/4 & \pm 3\sqrt{3}/4 & 9/4   \\
                                     \hline  4 & k & \pm\sqrt{k(3-k)} & -k+2 & \mp\sqrt{k(3-k)} & k   \\
                                     \hline  5  &-3k^2\pm\frac{3}{4}k+\frac{9}{4} & \frac{-3\sqrt{3}}{4}(3k\mp 1) & -k^2\pm\frac{1}{4}k-\frac{1}{4} &\frac{-3\sqrt{3}}{4}(k\pm 1) &-3k^2\mp\frac{21}{4}k+\frac{9}{4}  \\
                                     \hline  6 & k & \pm\sqrt{3}/2(2k-3) & 3k-1 & \pm 3\sqrt{3}/2(2k-1) & 9k-6   \\ \hline
                                    \end{array}$

\noindent In case 4, we have $0<k<3$.

%\noindent Note that the case $4$ can be regarded as the limit of case $6$ as $k$ tends to infinity.

As the proof of this theorem involves quite tedious computation, we write it as an appendix.

To classify the surfaces with constant curvatures, we need to consider the action of $D_3$ on $(I_{40},I_{31},I_{22},I_{13},I_{04})$.
\begin{lemma} The induced action of $D_3$ on $(I_{40},I_{31},I_{22},I_{13},I_{04})$ given by $\sigma$ is

$$\left(
\begin{array}{ccccc}
 I'_{40},I'_{31},I'_{22},I'_{13},I'_{04}
     \end{array}
   \right)=\frac{1}{16}\left(
\begin{array}{ccccc}
 I_{40},I_{31},I_{22},I_{13},I_{04}
   \end{array}
   \right)\left(
             \begin{array}{ccccc}
              1 & \sqrt{3} & 3 & 3\sqrt{3} & 9\\
              -4\sqrt{3} & -8 & -4\sqrt{3} & 0 & 12\sqrt{3}\\
              18 & 6\sqrt{3} & -2 & -6\sqrt{3} & 18\\
              -12\sqrt{3} & 0& 4\sqrt{3} & -8 & 4\sqrt{3}\\
              9 & -3\sqrt{3} & 3 & -\sqrt{3} & 1\\
             \end{array}
           \right).$$
%$I'_{40}=\frac{1}{16}(-12I_{13}\sqrt{3}-4I_{31}\sqrt{3}+I_{40}+9I_{04}+18I_{22}), I'_{31}=\frac{1}{16}(-3I_{04}\sqrt{3}+6I_{22}\sqrt{3}+I_{40}\sqrt{3}-8I_{31})$,

%$I'_{22}=\frac{1}{16}(4I_{13}\sqrt{3}-4I_{31}\sqrt{3}-2I_{22}+3I_{04}+3I_{40}), I'_{13}=\frac{1}{16}(-6I_{22}\sqrt{3}+3I_{40}\sqrt{3}-I_{04}\sqrt{3}-8I_{13})$,

%$I'_{04}=\frac{1}{16}(12I_{31}\sqrt{3}+4I_{13}\sqrt{3}+I_{04}+9I_{40}+18I_{22})$.

\noindent The induced action given by $\tau$ is
$$\left(
\begin{array}{ccccc}
 I'_{40},I'_{31},I'_{22},I'_{13},I'_{04}
     \end{array}
   \right)=\left(
\begin{array}{ccccc}
 I_{40},I_{31},I_{22},I_{13},I_{04}
   \end{array}
   \right)\left(
             \begin{array}{ccccc}
              1 & 0 & 0 & 0 & 0\\
              0 & -1 & 0 & 0 & 0\\
              0 & 0 & 1 & 0 & 0\\
              0 & 0& 0 & -1 & 0\\
              0 & 0 & 0 & 0 & 1\\
             \end{array}
           \right).$$

%$(I'_{40},I'_{31},I'_{22},I'_{13},I'_{04})=(I_{40},-I_{31},I_{22},-I_{13},I_{04})$.
\end{lemma}

\noindent{\bf Proof: }This is computed from Equation 17.

\begin{theorem}
Let $S$ be a regular elliptical surface with constant curvatures, then up to $D_3$ action, the curvature $(I_{40},I_{31},I_{22},I_{13},I_{04})$ of $S$ has one of the following forms.

1. $A^{1}_k =(k,0,-1,0,3)$, $A^2_k=(-3k^2+6k,0,-k^2-2k+2,0,-3k^2-6k)$ or $A^3_k=(k,0,k,0,-3k)$.

2. $B_{k,l}=(k,l(2k-3),k,0,-3k)$, where $l(2k-3)> 0$.

3. $C_{k}=(k,\sqrt{k(3-k)},-k+2,-\sqrt{k(3-k)},k)$.

%4. $D=(\frac{9}{4},\frac{3\sqrt{3}}{4},\frac{-1}{4} , \frac{- 3\sqrt{3}}{4},\frac{9}{4}  )$.
\end{theorem}

\noindent {\bf Proof}: This is computed from the previous theorem by considering the $D_3$ action. The $A^{1}_k$ form comes from the first case of Theorem 5.1, the $A^{2}_k$ form comes from the second case and the $A^{3}_k$ form comes from the third case with $l=0$. The $B_{k,l}$ form comes from the third case of Theorem 5.1 with $l\neq 0$. The $C_{k}$ form comes from the fourth case of Theorem 5.1. Note that under the action of $D_3$, case 1 and 6 in Theorem 5.1 can be transformed to each other by the action of $\sigma$ in Lemma 5.1. Case 2 and case 5 can be transformed to each other by the action of $\sigma$ in Lemma 5.1.

%The action of $\sigma$ on $B_{k,l}$ can be computed by the following result.

We need the following lemma to finish the classification of elliptical surfaces with constant curvatures.

\begin{lemma}
If $w_1,w_2$ are 1-forms satisfying the equation \begin{equation}dw_1=\lambda w_1\wedge w_2, dw_2=\mu w_1\wedge w_2, w_1\wedge w_2\not=0.\end{equation} Where $(\lambda,\mu)\in \mathbb{R}^2$ are constants. Then $w_1$ and $w_2$ are determined by $\lambda,\mu$ up to local diffeomorphisms.
\end{lemma}

\noindent {\bf Proof}: Case 1. If $(\lambda,\mu)=(0,0)$, the equations have the local solutions $w_1=df,w_2=dg$ with $df\wedge dg\not=0$, and $w_1,w_2$ are the pull back of $du,dv$ by the local diffeomorphism $u=f(x,y),v=g(x,y)$.

Case 2. If $(\lambda,\mu)=(\lambda,0)$, $\lambda\not=0$, the solutions of these equations must be of the forms $w_1=a df, w_2=dg$ for certain $f,g$ and $a$ with $adf\wedge dg\not=0$. Let $u=f(x,y),v=g(x,y)$, we get $-\frac{\partial{a}}{\partial v}=\lambda a$.
Hence we have $w_1=e^{-\lambda v}du, w_2=dv$ up to local diffeomorphisms.

Case 3. If $\mu \not=0$, define $w'_1=\frac{1}{\mu}w_2, w'_2=-\mu w_1+\lambda w_2$, then it is easy to check that $dw'_1=w'_1\wedge w'_2, dw'_2=0,w'_1\wedge w'_2=w_1\wedge w_2\not=0$. So $w'_1, w'_2$ can be determined by Case 2. It follows that $w_1$ and $w_2$ are determined up to local diffeomorphisms.

%\begin{lemma}
%If $(\lambda,\mu)\not=(0,0)$, then there exist $\lambda', \mu'$ such that $w_1,w_2$ is the solution to the Equations 24 with $(\lambda,\mu)=(1,0)$ if and only if $w'_1=\lambda w_1+\lambda'w_2, w'_2=\mu w_1+\mu'w_2$ is the solution to Equations 24 with $(\lambda,\mu)$.
%\end{lemma}

%If $(\lambda,\mu)\not=(0,0)$, then there exist $\lambda', \mu'$ such that $\lambda\mu'-\lambda'\mu=1$, they satisfy the condition of the lemma.

%This discussion shows that we have

%\begin{lemma}
%If $w_1,w_2$ and $w'_1,w'_2$ are two solutions of Equation 25, then there exists $f:\mathbb{R}^2\to \mathbb{R}^2$ such that $w'_1=f^*(w_1),w'_2=f^*(w_2)$ locally.
%\end{lemma}

\begin{theorem}
A regular elliptical surface with constant curvatures is determined by the invariants $I_{40},I_{31},I_{22},I_{13},I_{04}$ up to affine congruence.
\end{theorem}

\noindent {\bf Proof: }Let $S$ be a regular elliptical surface with constant curvatures $I_{40},I_{31},I_{22},I_{13},I_{04}$, then $w_1,w_2$ satisfy Equation 26 with $(\lambda,\mu)=(Y_1,Y_2)$,
where $\D{Y_1=\frac{2I_{13}}{3},Y_2=\frac{I_{40}}{4}-\frac{I_{22}}{2}-\frac{I_{04}}{12}}$ are defined in Section 4. By Lemma 5.2, $w_1,w_2$ are determined up to local diffeomorphisms.
Since $w_1,w_2,I_{40},I_{31},I_{22},I_{13},I_{04}$ are complete differential invariants, the theorem is proved.
%\begin{lemma}
%The action of $\sigma$ on $B_{kl}$ is given by

%$k'=\frac{\sqrt{3}}{2}lk-\frac{3\sqrt{3}}{4}l-\frac{1}{2}k, l'=\frac{-2\sqrt{3}k-2lk+3l}{2lk\sqrt{3}-3\sqrt{3}l-2k-6}$.

%and $k'=k,l'=-l$. And the fundamental domain of this action is $ $.
%\end{lemma}

\section{Appendix: the solutions of compatible equations}

To prove Theorem 5.1, We need to solve the compatible equations in Corollary 4.1. Denote

$J_1:=( -23I_{40}-6I_{22}+I_{04}+48) I_{13}+ ( -3I_{40}+18I_{22}+5I_{04}) I_{31}=0$

$J_2:=( I_{40}+6I_{22}-3I_{04})I_{13}+ (-3I_{40}+6I_{22}+I_{04}) I_{31}=0$

$K_1:=-9I_{40}I_{22}+3I_{40}I_{04}-24I_{31}I_{13}+18I_{22}^2-3I_{22}I_{04}+8I_{13}^2-I_{04}^2-18I_{40}-36I_{22}-18I_{04}=0$

$K_2:=-3I_{40}I_{22}-7I_{40}I_{04}-16I_{31}I_{13}+18I_{22}^2+15I_{22}I_{04}-16I_{13}^2+I_{04}^2+18I_{40}-36I_{22}-6I_{04}=0.$

$H_{11}:= -23I_{40}-6I_{22}+I_{04}+48, \ \ \ \ H_{12}:= -3I_{40}+18I_{22}+5I_{04}$

$H_{21}:=I_{40}+6I_{22}-3I_{04}, \ \ \ \ \ \ \ \ \ \ \ H_{22}:=-3I_{40}+6I_{22}+I_{04}$

$K_3:=H_{11}H_{22}-H_{12}H_{21}$
%\frac{1}{8}[( -6I_{22}-23I_{40}+I_{04}+48)(I_{04}+6I_{22}-3I_{40})- ( -3I_{40}+5I_{04}+18I_{22})( -3I_{04}+6I_{22}+I_{40})]$

$\ \ \ \ \  = 9I_{40}^2-15I_{40}I_{22}-5I_{40}I_{04}-18I_{22}^2+3I_{22}I_{04}+2I_{04}^2-18I_{40}+36I_{22}+6I_{04}$.

$K_{21}:=K_2-K_1=(6 I_{40}+18 I_{04})I_{22}-10 I_{40}I_{04}+8 I_{31}I_{13}-24 I_{13}^{2} +2 I_{04}^{2}+36 I_{40}+12 I_{04}=0$.

Case 1. If $K_3\not=0$, then by Equations 22 and 23, we must have $I_{13}=I_{31}=0$.

By substituting $I_{13}=I_{31}=0$ in $K_{21}$, we obtain
\begin{equation}(3I_{40}+9I_{04})I_{22}-5I_{40}I_{04}+I_{04}^2+18I_{40}+6I_{04}=0\end{equation}

i) If $I_{40}+3I_{04}=0$, then $-5I_{40}I_{04}+I_{04}^2+18I_{40}+6I_{04}=16I^2_{04}-48I_{04}=0$. We have $(I_{40},I_{04})=(0,0)$ or $(-9,3)$.

If $(I_{40},I_{04})=(0,0)$, then by $K_1=0$ we get $I_{22}=0$ or $2$. In this case

$(I_{40},I_{31},I_{22},I_{13},I_{04})=(0,0,0,0,0)$ or $(0,0,2,0,0)$.

These two solutions are included in the third and the second cases of the table in Theorem 5.1 respectively.

If $(I_{40},I_{04})=(-9,3)$, then by $K_1=0$ we get $I_{22}=-1$. In this case

$(I_{40},I_{31},I_{22},I_{13},I_{04})=(-9,0,-1,0,3)$.

This solution is included in the first case of the table in Theorem 5.1.

ii) If $I_{40}+3I_{04}\not =0$, then by Equation 27, we have \begin{equation}I_{22}=-\frac{-5I_{40}I_{04}+I_{04}^2+18I_{40}+6I_{04}}{3(I_{40}+3I_{04})}\label{I_{22}}.\end{equation}

By substituting the expression of $I_{22}$ into $K_1$, we get

$$K_1=-\frac{4(I_{04}-3)(3I_{40}+I_{04})(I_{40}^2-2I_{40}I_{04}+I_{04}^2+24I_{40}+24I_{04})}{(I_{40}+3I_{04})^2}.$$
Hence $I_{04}=3$, $ 3I_{40}+I_{04}=0$ or $I_{40}^2-2I_{40}I_{04}+I_{04}^2+24I_{40}+24I_{04}=0$.

If $I_{04}=3$, then by Equation 28 we get $I_{22}=-1$. Set $I_{40}=k$,
%By $K_3\not=0$, we get $I_{40}^2-2I_{40}-3\not=0$.
 we have

 $(I_{40},I_{31},I_{22},I_{13},I_{04})=(k,0,-1,0,3)$.

 This solution is included in the first case of the table in Theorem 5.1.

If $3I_{40}+I_{04}=0$, then $I_{04}=-3I_{40}$. By Equation 28 we get $I_{22}=I_{40}$. Set $I_{40}=k$, we have

$(I_{40},I_{31},I_{22},I_{13},I_{04})=(k,0,k,0,-3k)$.

This solution is included in the third case of the table in Theorem 5.1.

If $I_{40}^2-2I_{40}I_{04}+I_{04}^2+24I_{40}+24I_{04}=0$, then by regarding this equation as a parabola, we can parameterize it as $I_{40}=-3k^2+6k, I_{04}=-3k^2-6k$. By Equation 28 we get $I_{22}=-k^2-2k+2$.
Hence

$(I_{40},I_{31},I_{22},I_{13},I_{04})=(-3k^2+6k,0,-k^2-2k+2,0,-3k^2-6k)$.

This solution is included in the second case of the table in Theorem 5.1.
%By substituting the expression of $I_{22}$ into $K_3$ we get $$K_3=\frac{(I_{04}+3I_{40})(13I_{04}^3-23I_{04}^2I_{40}+7I_{04}I_{40}^2+3I_{40}^3-24I_{04}^2+192I_{04}I_{40}+24I_{40}^2-288I_{04}-288I_{40})}{(3I_{04}+I_{40})^2}$$ and we need $K_3\not=0$.

Case 2. If $K_3=0$, then $I_{13}$ or $I_{31}$ can be nonzero.

i) If $H_{11}=-23I_{40}-6I_{22}+I_{04}+48=0$ and $H_{12}=-3I_{40}+18I_{22}+5I_{04}=0$, then $\D{I_{04}=9I_{40}-18, I_{22}=-\frac{7}{3}I_{40}+5}$. Equation 23 can be written as $H_{21}I_{13}+ H_{22}I_{31}=0$. As in this case it is easy to check that $H_{21}$ and $H_{22}$ can not be both zero. Then $I_{13}=kH_{22}=k(-3I_{40}+6I_{22}+I_{04}), I_{31}=-k H_{21}=-k( I_{40}+6I_{22}-3I_{04})$, by substituting the expression into $K_1$ and $K_2$ we get

$K_1=4(2I_{40}-3)(1024I_{40}k^2-2112k^2+16I_{40}-45)=0,$

$K_2=8(64k^2-3)(2I_{40}-3)(4I_{40}-9)=0.$

By solving these two equations we get solutions $(I_{40},k)=(\frac{3}{2},k),(\frac{9}{4},\pm\frac{{\sqrt3}}{8})$. Hence

$(I_{40},I_{31},$ $I_{22},I_{13},I_{04})=(\frac{3}{2},-24k,\frac{3}{2},0,-\frac{9}{2})$ or $(\frac{9}{4},\frac{\mp 3\sqrt{3}}{4},-\frac{1}{4},\frac{\pm 3\sqrt{3}}{4},\frac{9}{4})$.

These two solutions are included in the third and fourth cases of the table in Theorem 5.1 respectively.

ii) If at least one of $H_{11},H_{12}$ is not zero. Equation 22 can be written as $H_{11}I_{13}+ H_{12}I_{31}=0$. Then $I_{13}=kH_{12}=k( -3I_{40}+18I_{22}+5I_{04}), I_{31}=-kH_{11}=-k( -23I_{40}-6I_{22}+I_{04}+48)$. By substituting the expression into $K_1,K_2$ and $K_3$ we get

$K_1=k^2(1728I_{40}^2-10368I_{40}I_{22}-3072I_{40}I_{04}+1152I_{22}I_{04}+320I_{04}^2-3456I_{40}+20736I_{22}+5760I_{04})$
\noindent$-9I_{40}I_{22}+3I_{40}I_{04}+18I_{22}^2-3I_{22}I_{04}-I_{04}^2-18I_{40}-36I_{22}-18I_{04}=0$

$K_2=k^2(960I_{40}^2-4608I_{40}I_{22}-1408I_{40}I_{04}-6912I_{22}^2-3072I_{22}I_{04}-320I_{04}^2-2304I_{40}+13824I_{22} \noindent +3840I_{04})
-3I_{40}I_{22}-7I_{40}I_{04}+18I_{22}^2+15I_{22}I_{04}+I_{04}^2+18I_{40}-36I_{22}-6I_{04}=0$

$K_3=9I_{40}^2-15I_{40}I_{22}-5I_{40}I_{04}-18I_{22}^2+3I_{22}I_{04}+2I_{04}^2-18I_{40}+36I_{22}+6I_{04}$.

Rewrite $K_1$ and $K_2$ by $K_1=k^2 G_{11}+G_{12}$ and $K_2=k^2 G_{21}+G_{22}$, where $G_{11},G_{12},G_{21},G_{22}$ are polynomials in $I_{40},I_{22},I_{04}$. We set $L=G_{11}G_{22}-G_{12}G_{21}$. Then $K_1=0$ and $K_2=0$ make $L=0$. Computation shows

$L=384(-3I_{40}+18I_{22}+5I_{04})(-3I_{40}^2I_{22}+13I_{40}^2I_{04}-21I_{40}I_{22}^2-24I_{40}I_{22}I_{04}-3I_{40}I_{04}^2+18I_{22}^3+3I_{22}^2I_{04}+I_{22}I_{04}^2
-42I_{40}^2+15I_{40}I_{22}-33I_{40}I_{04}-18I_{22}^2+21I_{22}I_{04}+I_{04}^2+90I_{40}-36I_{22}+18I_{04})=0$

Set $M=-3I_{40}^2I_{22}+13I_{40}^2I_{04}-2I_{40}1I_{22}^2-24I_{40}I_{22}I_{04}-3I_{40}I_{04}^2+18I_{22}^3+3I_{22}^2I_{04}+I_{22}I_{04}^2-42I_{40}^2+15I_{40}I_{22}-33I_{40}I_{04}
-18I_{22}^2+21I_{22}I_{04}+I_{04}^2+90I_{40}-36I_{22}+18I_{04}$. Then $L=384H_{12}M$.

a) If $H_{12}=-3I_{40}+18I_{22}+5I_{04}=0$, then $H_{11}=-23I_{40}-6I_{22}+I_{04}+48\not=0$. Since $J_1=0$ we have $I_{13}=0$. And $\D{I_{22}= \frac{3I_{40}-5I_{04}}{18}}$;
By substituting $I_{22}$ into $K_1,K_2$ and $K_3$, we have

$K_1=\frac{1}{9}(3I_{40}+I_{04})(-3I_{40}+11I_{04}-72)=0$

$K_2=-\frac{4}{9}(3I_{40}+I_{04})(4I_{04}-9)=0$

$K_3=-\frac{2}{9}(3I_{40}+I_{04})(-9I_{40}+I_{04}+18)=0$

Since $K_1=K_2=K_3=0$, we must have $3I_{40}+I_{04}=0$. And we have a solution of the form

$(I_{40},I_{31},I_{22},I_{13},I_{04})=(I_{40},l(2I_{40}-3),I_{40},0,-3I_{40})$.

This solution is included in the third case of the table in Theorem 5.1.

b) If $M=0$, then by $K_3=0$ we have

$I_{04}^2=-\frac{9}{2} I_{40}^{2}+\frac{15}{2}I_{40}I_{22}+\frac{5}{2}I_{40}I_{04}+9 I_{22}^{2}-\frac{3}{2}I_{22}I_{04}+9I_{40}-18 I_{22}-3I_{04}$.

By substituting the expression of $I_{04}^2$ into $M=0$, we have

$M=\frac{1}{2}(-I_{40}+3I_{22}+3)(-27I_{40}^2-21I_{40}I_{22}-11I_{40}I_{04}+18I_{22}^2+I_{22}I_{04}+66I_{40}-36I_{22}+10I_{04})=0$

If $-I_{40}+3I_{22}+3=0$, then $I_{40} =3I_{22}+3$. By substituting $I_{40}$ into $K_3$, we get $K_3=18I_{22}^2-12I_{22}I_{04}+2I_{04}^2+99I_{22}-9I_{04}+27 =0$. Regarding this equation as a parabola equation, it can be  parameterized as $$\D{I_{04}=-\frac{3}{8}(2l^2-11l+3), I_{22}=-\frac{1}{8}(-2l^2-3l+3)}.$$ Inserting the expression of $I_{40},I_{22},I_{04}$ into $K_1,K_2$, we have
$K_1=-9l (l-3) (256(l-1)^2 k^2-3)$ and $K_2=-9(l-1)(l-3)(256(l-1)^2 k^2-3)$. Hence $l=3$ or $256(l-1)^2 k^2=3$. And we have

$(I_{40},I_{31},I_{22},I_{13},I_{04})=(-3/2,-96k,-3/2,0,9/2)$ or
$\D{(\frac{3(768k^2\pm8\sqrt{3}k-3)}{1024k^2},\frac{3\sqrt{3}(\pm32k-3\sqrt{3})}{128k}}$, $\D{\frac{(-256k^2\pm8\sqrt{3}k-3)}{1024k^2},
\frac{3\sqrt{3}(\mp32k-\sqrt{3})}{128k},\frac{3(768k^2\mp56\sqrt{3}k-3)}{1024k^2})}$.

We replace $k$ by $\frac{\sqrt{3}}{32k}$, then the later case can be written as

$(I_{40},I_{31},I_{22},I_{13},I_{04})=(-3k^2\pm\frac{3}{4}k+\frac{9}{4},\frac{-3\sqrt{3}}{4}(3k\mp 1),-k^2\pm\frac{1}{4}k-\frac{1}{4},\frac{-3\sqrt{3}}{4}(k\pm 1),-3k^2\mp\frac{21}{4}k+\frac{9}{4})$.

 These solutions are included in the third and the fifth cases of the table in Theorem 5.1 respectively.

%$K_1=-81l(l-1)(256(3l-1)^2 k^2-3)$ $K_2=-27(3l-1)(l-1)(256(3l-1)^2 k^2-3)$

%$K_1=(320I_{04}^2-8064I_{04}I_{22}-15552I_{22}^2-3456I_{04}+10368I_{22}+5184)k^2-I_{04}^2+6I_{04}I_{22}-9I_{22}^2-9I_{04}-117I_{22}-54=0$

%$K_2=(-320I_{04}^2-7296I_{04}I_{22}-12096I_{22}^2-384I_{04}+10368I_{22}+1728)k^2+I_{04}^2-6I_{04}I_{22}+9I_{22}^2-27I_{04}+9I_{22}+54=0$

%By $K_1+K_2=0$, we solve $K=k^2=-3/64(I_{04}+3I_{22})/(4I_{22}+1)/(-9+9I_{22}+5I_{04})$. Substituting $k^2$ into $K_1$ and $K_2$ we can show

%$\frac{-3(3+8\sqrt{3}k-768k^2)}{1024k^2},\frac{-3\sqrt{3}(32k+3\sqrt{3})}{128k},\frac{-(3+8\sqrt{3}k+256k^2)}{1024k^2},\frac{-3\sqrt{3}(-32k+\sqrt{3})}{128k},\frac{3(56\sqrt{3}k+768k^2-3)}{1024k^2}$

%$\frac{3(768k^2\pm8\sqrt{3}k-3)}{1024k^2},\frac{-3\sqrt{3}(-32k+3\sqrt{3})}{128k},\frac{(-3+8\sqrt{3}k-256k^2)}{1024k^2},\frac{-3\sqrt{3}(32k+\sqrt{3})}{128k},\frac{-3(3+56\sqrt{3}k-768k^2)}{1024k^2}$

%$(I_{40},\pm \sqrt{K}(-75I_{22}-21+I_{04}),I_{22},-\sqrt{K}(-9+9I_{22}+5I_{04}),I_{40})$ with $2I_{04}^2-12I_{04}I_{22}+18I_{22}^2-9I_{04}+99I_{22}+27 =0$ is a solution.

If $-27I_{40}^2-21I_{40}I_{22}-11I_{40}I_{04}+18I_{22}^2+I_{22}I_{04}+66I_{40}-36I_{22}+10I_{04}=0$, then by adding
$K_3=9I_{40}^2-15I_{40}I_{22}-5I_{40}I_{04}-18I_{22}^2+3I_{22}I_{04}+2I_{04}^2-18I_{40}+36I_{22}+6I_{04}=0$ to both side of the above equation, we have
$(-36I_{40}+4I_{04})I_{22}-18I_{40}^2-16I_{40}I_{04}+2I_{04}^2+48I_{40}+16I_{04}=0$.

If $-36I_{40}+4I_{04}=0$, then we get

$(I_{40},I_{31},I_{22},I_{13},I_{04})=(0,-48k,0,0,0)$.

This solution is included in the third case of the table in Theorem 5.1.

If $-36I_{40}+4I_{04}\not=0$, then we have
$$\D{I_{22}=-\frac{-9I_{40}^2-8I_{40}I_{04}+I_{04}^2+24I_{40}+8I_{04}}{-18I_{40}+2I_{04}}}.$$ By substituting $I_{22}$ into $K_3=0$, we have
$$\D{K_3=-\frac{4(-I_{40}+I_{04})(3I_{40}+I_{04})(-9I_{40}+I_{04}+6)(-9I_{40}+I_{04}+18)}{(-9I_{40}+I_{04})^2}}=0.$$

If $I_{40}-I_{04}=0$, then we get

$(I_{40},I_{31},I_{22},I_{13},I_{04})=(I_{40},4k(4I_{40}-9),-I_{40}+2,-4k(4I_{40}-9),I_{40})$.

By $K_1=0$, we have $k^2=-\frac{1}{16}I_{40}(I_{40}-3)/(4I_{40}-9)^2$.
This solution is the fourth case of the table in Theorem 5.1.

If $3I_{40}+I_{04}=0$, then we get

$(I_{40},I_{31},I_{22},I_{13},I_{04})=(I_{40},16k(2I_{40}-3),I_{40},0,-3I_{40})$.

This solution is included in the third case of the table in Theorem 5.1.

If $-9I_{40}+6+I_{04}=0$, then we get

$(I_{40},I_{31},I_{22},I_{13},I_{04})=(\frac{1}{2},-32k ,\frac{1}{2},0,-\frac{3}{2})$ or $(I_{40},\pm\sqrt{3}/2(2I_{40}-3),3I_{40}-1,\pm 3\sqrt{3}/2(2I_{40}-1),9I_{40}-6)$.

These solutions are included in the third and the sixth cases of the table in Theorem 5.1.

If $-9I_{40}+18+I_{04}=0$, then we get

$(I_{40},I_{31},I_{22},I_{13},I_{04})=(\frac{3}{2},0,\frac{3}{2},0,-\frac{9}{2})$.

This solution is included in the third case of the table in Theorem 5.1.

These are all possible cases for elliptical surfaces with constant curvatures.

\end{document}